\documentclass[a4paper, fullpage, reqno, 11pt]{amsart}%,draft
\usepackage{amssymb}
\usepackage{ifthen} % pachetto necessario per un minimo di
\provideboolean{shownotes} % define a boolean variable
\setboolean{shownotes}{true} % defined as ``true'' or ``false''
\newcommand{\margnote}[1]{
\ifthenelse{\boolean{shownotes}}%
{\marginpar{\raggedright\tiny\texttt{#1}}}%
{}%
}

\newcommand{\hole}[1]{
\ifthenelse{\boolean{shownotes}}%
{\begin{center} \fbox{ \rule {.25cm}{0cm}
\rule[-.1cm]{0cm}{.4cm} \parbox{.85\textwidth}{\begin{center}
\texttt{#1}\end{center}} \rule {.25cm}{0cm}}\end{center}}
{}
}
\newtheorem{theorem}{Theorem}[section]
\newtheorem{proposition}[theorem]{Proposition}
\newtheorem{lemma}[theorem]{Lemma}
\newtheorem{corollary}[theorem]{Corollary}

\theoremstyle{remark}

\newtheorem{definition}[theorem]{Definition}

\newcommand{\e}{\varepsilon}		       
\newcommand{\R}{\mathbb{R}}

\newcommand{\ue}{u^{\varepsilon}}
\newcommand{\ut}{\tilde{u}}
\newcommand{\pe}{p^{\varepsilon}}
\newcommand{\pt}{\tilde{p}}
\newcommand{\te}{\theta^{\varepsilon}}
\newcommand{\tht}{\tilde{\theta}}
\newcommand{\dive}{\mathop{\mathrm {div}}}
\numberwithin{equation}{section}

\begin{document}
\title [Artificial compressibility for Navier Stokes Fourier]
{On the artificial compressibility method for the Navier Stokes Fourier system}
%{A DISPERSIVE APPROACH \\TO THE ARTIFICIAL COMPRESSIBILITY APPROXIMATIONS \\OF THE NAVIER STOKES EQUATIONS IN 3-D}

\author{Donatella Donatelli}
\address{Donatella Donatelli --- 
Dipartimento di Matematica Pura ed Applicata \\
Universit\`a di L'Aquila \\
Via Vetoio, \\
		     67010  Coppito (AQ), Italy}
\email{donatell@univaq.it}

%\author{Pierangelo Marcati}
%\address{Pierangelo Marcati --- 
%Dipartimento di Matematica Pura ed Applicata \\
%Universit\`a di L'Aquila \\
%Via Vetoio, \\
%		     67010  Coppito (AQ), Italy}
%\email{marcati@univaq.it}

\begin{abstract}
This paper  deals with the  approximation of the  weak solutions of the incompressible Navier Stokes Fourier system. In particular it extends the artificial compressibility method  for the Leray weak solutions of the Navier Stokes equation, used by Temam \cite{Tem01},  in the case of a bounded domain and later in \cite{DM06} in the case of the whole space. By exploiting the wave equation structure of the pressure  of the approximating system  the convergence of the approximating sequences is achieved  by means of dispersive  estimate of Strichartz type. It will be proved that the projection of the approximating  velocity fields on the divergence free vectors is relatively compact and converges to a  weak solution of the incompressible Navier Stokes Fourier system. \end{abstract}
\date{}
\maketitle
% \baselineskip 18pt
% \tableofcontents
\section{Introduction}\label{intro}

This paper is concerned with  the  artificial compressibility approximation for the following Navier Stokes Fourier system, namely
\begin{equation}
    \begin{cases}
	\partial_{t}u+(u\cdot\nabla)u-\mu\Delta u=\nabla p\\
	\partial_{t}\theta+u\cdot\nabla\theta=\kappa\Delta \theta\\
	\dive u=0\\
	u(x,0)=u_{0}(x),\ \theta(x,0)=\theta_{0}(x)
	\end{cases}
    \label{1}
\end{equation}
where $(x,t)\in \R^{3}\times [0,T]$, $u\in\R^{3}$ denotes the velocity vector field , $p\in \R$ the pressure of the fluid , $\theta\in\R$ is the temperature of the fluid, $\mu$ is the kinematic viscosity and $\kappa$ is the heat conductivity. The model \eqref{1} was proposed in \cite{GSR04} and \cite{LPL96} and, as we can observe the system \eqref{1} is an extension of the Navier Stokes equation which governs both the velocity fluid  $u$ and the fluctuation of the temperature field $\theta$ in an incompressible fluid. Let us point out that in the model under consideration the temperature field is just advected by the velocity field $u$ and diffuses according to Fourier's law. More complicated effects such as viscous heating for example (see \cite{LPL96}) will be not considered in the present paper, but they will be treated in future work. The artificial compressibility approximation was introduced by Chorin \cite{Ch68, Ch69}, Temam \cite{Tem69a, Tem69b} and Oskolkov \cite{Osk}, in order to deal with the difficulty induced by the incompressibility constraints in the numerical approximations to the Navier Stokes equation. The paper of Temam \cite{Tem69a, Tem69b} and his book \cite{Tem01} discuss the convergence of these approximations on bounded domains by using the classical Sobolev compactness embedding theorems and they recover the compactness in time by the classical Lions \cite{L-JL59} method of fractional derivatives. In \cite{DM06} the authors cannot make use of the classical compactness theorems, therefore they exploit the underlying wave equation structure and recover the necessary compactness by means  of dispersive type estimates. Here it will be used   a similar approach adapted to the artificial compressibility approximation of the Navier Stokes Fourier system. In particular the following system we will considered
\begin{equation}
\begin{cases}
\displaystyle{\partial_{t}\ue+\nabla \pe=\mu\Delta \ue-\left(\ue\cdot\nabla\right)\ue -\frac{1}{2}(\dive \ue)\ue}\\
\displaystyle{\partial_{t}\te+\ue\cdot\nabla\te=\kappa\Delta \te -\frac{1}{2}(\dive \ue)\te}\\
\e\partial_{t}\pe+ \dive \ue=0,
\end{cases}
\label{2}
\end{equation}
where $(x,t)\in \R^{3}\times [0,T]$, $\ue=\ue(x,t)\in\R^{3}$,   $\te=\te(x,t)\in\R$ and $\pe=\pe(x,t)\in \R$.\\
The system will be regarded as a semilinear wave type equation for the pressure function and the dispersive estimates will be carried out by using certain classical  $L^{p}$-type estimates due to Strichartz \cite{GV95, KT98, S77}. The particular type of  Strichartz estimates that we are going to use here can be found in the book of Sogge \cite{So95} or deduced by the so called bilinear estimates of Klainerman and Machedon \cite{KM93} and  Foschi Klainerman \cite{FK00}. Our analysis can also be related to the convergence of the incompressible limit problem via a formal expansion (see for instance Temam \cite{Tem01}, Chapter 3). In particular our  wave equation structure has some similarities with the one  exploited in various way by the paper of  Desjardin, Grenier, Lions, Masmoudi \cite{DGLM}, Desjardin and Grenier \cite{DG99}.\\
The interest into the artificial compressibility methods started with the previously mentioned results of Chorin and Temam and was later on investigated by Ghidaglia and Temam \cite{GhT88}. Later developments of numerical investigations in the directions of projections methods have been carried out by \cite{GMS06}, \cite{E03}, \cite{Pro97}, \cite{Ran92}, \cite{NP}, \cite{KS}, \cite{Sma}, \cite{YKS}.\\
The approximation of the Navier Stokes Fourier system has also been analyzed by a different point of view. In fact F. Golse and L. Saint-Raimond in \cite{GSR04}, \cite{GSR05} recover weak solutions for the system \eqref{1} in the frame of the hydrodynamic limit for the Boltzmann equation.\\
This paper is organized as follows. In Section 2 we recall the mathematical tools needed in the paper and we recall same basic definitions. In Section 3 we set up our problem, we explain our approximating system and we state our main result. The Section 4 is devoted to recover the a priori estimates needed to  get the strong convergence of the approximating sequences and to prove the main theorem. Finally, in Section 5, we give the proof of the main result.

\section{Preliminaries}

For convenience of the reader we establish some notations and recall some basic facts that will be useful in the sequel.\\
We will denote by  $\mathcal{D}(\R ^d \times \R_+)$
 the space of test function
$C^{\infty}_{0}(\R^d \times \R_+)$, by $\mathcal{D}'(\R^d \times
\R_+)$ the space of Schwartz distributions and $\langle \cdot, \cdot \rangle$
the duality bracket between $\mathcal{D}'$ and $\mathcal{D}$ and by $\mathcal{M}_{t}X'$ the space $C_{c}^{0}([0,T];X)'$. Moreover
$W^{k,p}(\R^{d})=(I-\Delta)^{-\frac{k}{2}}L^{p}(\R^{d})$ and $H^{k}(\R^{d})=W^{k,2}(\R^{d})$ denote the nonhomogeneous Sobolev spaces for any $1\leq p\leq \infty$ and $k\in \R$. $\dot W^{k,p}(\R^{d})=(-\Delta)^{-\frac{k}{2}}L^{p}(\R^{d})$ and $\dot H^{k}(\R^{d})=W^{k,2}(\R^{d})$  denote the homogeneous Sobolev spaces. The notations
 $L^{p}_{t}L^{q}_{x}$ and $L^{p}_{t}W^{k,q}_{x}$ will abbreviate respectively  the spaces $L^{p}([0,T];L^{q}(\R^{d}))$, and $L^{p}([0,T];W^{k,q}(\R^{d}))$.\\
The operators $Q$ and $P$ denote   the Leray's projectors  on the space of gradients vector fields and on the space of divergence - free vector fields, respectively. Hence, in the sense of distribution, one has
\begin{equation}
Q=\nabla \Delta^{-1}\dive\qquad P=I-Q.
\label{3}
\end{equation} 
Let us remark that   $Q$ and $P$  can be expressed in terms of Riesz multipliers, therefore they are  bounded linear operators on every $W^{k,p}$ $(1<p<\infty)$ space (see \cite{Ste93}).   \\ \\
Let us recall that if  $w$ is a (weak) solution of the following wave equation in the space $[0,T]\times \R^{d}$
\begin{equation*}
\begin{cases}
\left(-\frac{\partial ^{2}}{\partial t}+\Delta\right)w(t,x)=F(t,x)\\
w(0,\cdot)=f,\quad \partial_{t}w(0,\cdot)=g,
\end{cases}
\end{equation*}
for some data $f,g, F$ and time $0<T<\infty$, 
then $w$ satifies the following Strichartz estimates, (see \cite{GV95}, \cite{KT98})
\begin{equation}
\|w\|_{L^{q}_{t}L^{r}_{x}}+\|\partial_{t}w\|_{L^{q}_{t}W^{-1,r}_{x}}\lesssim \|f\|_{\dot H^{\gamma}_{x}}+\|g\|_{\dot H^{\gamma -1}_{x}}+\|F\|_{L^{\tilde{q}'}_{t}L^{\tilde{r}'}_{x}},
\label{s2}
\end{equation}
where $(q,r)$, $(\tilde{q},\tilde{r})$ are \emph{wave admissible} pairs, namely they satisfy 
\begin{equation*}
\frac{2}{q}\leq (d-1)\left(\frac{1}{2}-\frac{1}{r}\right) \qquad 
\frac{2}{\tilde{q}}\leq (d-1)\left(\frac{1}{2}-\frac{1}{\tilde{r}}\right)
\end{equation*}
and moreover the following   conditions holds
\begin{equation*}
\frac{1}{q}+\frac{d}{r}=\frac{d}{2}-\gamma=\frac{1}{\tilde{q}'}+\frac{d}{\tilde{r}'}-2.
\end{equation*}
Later on we shall use \eqref{s2} in the case of $d=3$, $(\tilde{q}', {\tilde{r}'})=(1, 3/2)$,  so one has $\gamma=1/2$ and $(q,r)=(4,4)$, namely the following estimate holds
\begin{equation}
\|w\|_{L^{4}_{t,x}}+\|\partial_{t}w\|_{L^{4}_{t}W^{-1,4}_{x}}\lesssim \|f\|_{\dot H^{1/2}_{x}}+\|g\|_{\dot H^{ -1/2}_{x}}+\|F\|_{L^{1}_{t}L^{3/2}_{x}}.
\label{s3}
\end{equation}
Beside the Strichartz estimate \eqref{s2} or \eqref{s3} in the case of $d=3$ (see \cite{So95}), there exists a non standard estimate which follows from an earlier linear Strichartz \cite{S77} estimate. This inequality can also be deduced by the bilinear estimates of Klainerman and Machedon \cite{KM93}, Foschi  and Klainerman \cite{FK00}, namely
\begin{equation}
\|w\|_{L^{4}_{t,x}}+\|\partial_{t}w\|_{L^{4}_{t}W^{-1,4}_{x}}\lesssim \|f\|_{\dot H^{1/2}_{x}}+\|g\|_{\dot H^{1/2}_{x}}+\|F\|_{L^{1}_{t}L^{2}_{x}}.
\label{s1}
\end{equation}
Finally, we mention here the following lemma that will be usefull in the next sections,
\begin{lemma}
\label{ly}
Let us consider  a smoothing kernel $j\in C^{\infty}_{0}(\R^{d})$, such that $j\geq 0$, $\int_{\R^{d}}j dx=1$, and define the Friedrichs mollifiers as
\begin{equation*}
j_{\alpha}(x)=\alpha^{-d}j\left(\frac{x}{\alpha}\right).
\end{equation*}
Then  for any $f\in \dot H^{1}(\R^{d})$, one has
\begin{equation}
\label{y1}
\|f-f\ast j_{\alpha}\|_{L^{p}(\R^{d})}\leq C_{p}\alpha^{1-\sigma}\|\nabla f\|_{L^{2}(\R^{d})},
\end{equation}
where
\begin{equation*}
p\in [2, \infty)
\quad \text{if $d=2$}, \quad p\in [2, 6] \quad \text{if $d=3$ \ and}\quad \sigma=d\left(\frac{1}{2}-\frac{1}{p}\right).
\end{equation*}
Moreover the following Young type inequality holds
\begin{equation}
\label{y2}
\|f\ast j_{\alpha}\|_{L^{p}(\R^{d})}\leq C\alpha^{-s-d\left(\frac{1}{q}-\frac{1}{p}\right)}\|f\|_{W^{-s,q}(\R^{d})},
\end{equation}
for any $p,q\in [1, \infty]$, $q\leq p$,  $s\geq 0$, $\alpha\in(0,1)$.
\end{lemma}

\section{Approximating system and main result}
Let us consider the incompressible Navier Stokes Fourier system \eqref{1} we  recall (see P.L.Lions \cite{LPL96} ) the notion of  weak solution that will be used later for \eqref{1}.
\begin{definition}
\label{df}
We say that the couple $(u,\theta)$ is a weak solution of the Navier Stokes Fourier system if it  satisfies  
\begin{align*}
&\int_{0}^{T}\!\!\int_{\R^{d}}\left(\mu\nabla u\cdot\nabla\varphi -u_{i}u_{j}\partial_{i}\varphi_{j}-u\cdot
\frac{\partial \varphi}{\partial t}\right) dxdt
=\int_{\R^{d}}u_{0}\cdot \varphi dx,
\end{align*}
and
\begin{align*}
&\int_{0}^{T}\!\!\int_{\R^{d}}\left(\kappa\nabla \theta\cdot\nabla\varphi -\theta u\cdot\nabla\varphi-\theta\cdot
\frac{\partial \varphi}{\partial t}\right) dxdt
=\int_{\R^{d}}\theta_{0}\cdot \varphi dx,
\end{align*}
for all $\varphi\in C^{\infty}_{0}(\R^{d}\times[0,T])$, $\dive \varphi =0$ and
\begin{equation*}
\dive u=0 \qquad \text{in $\mathcal{D}'(\R^{d}\times[0,T])$}
\end{equation*}
Moreover the following energy inequality holds 
\begin{align}
\frac{1}{2}\int_{\R^{d}}&(|u(x,t)|^{2}+|\theta(x,t)|^{2})dx+\int_{0}^{t}\!\!\int_{\R^{d}}(\mu|\nabla u(x,s)|^{2}+\kappa|\nabla \theta(x,s)|^{2})dxds\notag\\ \leq
&\frac{1}{2}\int_{\R^{d}}(|u_{0}|^{2}+|\theta_{0}|^{2})dx, \qquad \text{for all $t\geq 0$}.
\label{4}
\end{align}
\end{definition}
Sometimes one refers also to this solution as Leray weak solutions although Leray himself didn't study thermal effects. It is worth to mention that the conservation of energy \eqref{4} is strictly related to the space dimension $d=3$ because of the integrability requirements of $|u|^{3}$ or $u\theta$. For a detailed mathematical description of  models involving temperature and on the connected mathematical difficulties in studying them we refer to the book of P.L.Lions \cite{LPL96}.
%There exists in the mathematical literature several results concerning the existence of Leray weak solutions to the Navier Stokes equations, for example we can refer to books of P.L.Lions \cite{LPL96} and Temam \cite{Tem01}. The case $d=3$ is a major open problem and a considerably more difficult case than the case $d=2$, since the bound on the $L^{2}$ norm (kinetic energy) provides only a control on a supercritical norm and does not provide any information concerning the critical controlling (and scaling 
%invariant) norm $L^{3}$. Hence we do not know (opposite to the case $d=2$) whether or not the Leray weak solutions are unique, unless (see Serrin \cite{Se63}) we assume a control on the $L^{3}$
%norm. Some important regularity results can be found in \cite{CKN}.\\
In order to approximate the system\eqref{1} we wish to use the system  \eqref{2}. As we can observe the first two equations of the system \eqref{2} resemble respectively the equation for conservation of momentum and the equation of the temperature balance for the system \eqref{1} except for the two terms  
\begin{equation*}
-\frac{1}{2}(\dive\ue)\ue, \qquad -\frac{1}{2}(\dive\ue)\te.
\end{equation*}
We need to add these  two correction in order to avoid the paradox of increasing the kinetic energy along the motion. In the third equation of the system \eqref{2} we introduce  a ``linearized'' compressibility constraint given by the equation
\begin{equation*}
\partial_{t}\pe=-\frac{1}{\e}\dive\ue.
\end{equation*}
The limiting behaviour as $\e\downarrow 0$ of the initial data of \eqref{2} deserves a little discussion. Indeed to solve \eqref{2},  three initial conditions are required
\begin{equation}
\ue(x,0)=\ue_{0}(x), \qquad \te(x,0)=\te_{0}(x),\qquad \pe(x,0)=\pe_{0}(x),
\end{equation}
while the Navier Stokes Fourier system \eqref{1} requires only an initial condition on the velocity $u$ and an initial condition on the temperature $\theta$. Hence our approximation will be consistent if the initial datum on the pressure will be eliminated by an ``initial layer'' phenomenon. Since in the limit we have to deal with weak solutions that verify \eqref{4} it is reasonable to require the finite energy constraint to be satisfied by the approximating sequences $(\ue,\te, \pe)$. So we can deduce a natural behaviour to be imposed on the initial data $(\ue_{0}, \te_{0}, \pe_{0})$, namely
\begin{align}
\tag {\bf{ID}}
 \ue_{0}&=\ue(\cdot, 0)\longrightarrow u_{0}=u(\cdot ,0)\ \text{strongly in}\  L^{2}(\R^{3})
\\ 
 \te_{0}&=\ue(\cdot, 0)\longrightarrow \theta_{0}=\theta(\cdot ,0)\ \text{strongly in}\  L^{2}(\R^{3})\notag\\
\sqrt{\e}\pe_{0}&=\sqrt\e \pe(\cdot,0)\longrightarrow 0\  \text{strongly in} \  L^{2}(\R^{3})\notag.
\end{align}
Let us remark that  the convergence of $\sqrt{\e}\pe_{0}$ to $0$ is necessary to avoid  the presence of concentrations of energy in the limit.\\
Since it will not affect our approximation process, for semplicity from now on, we will take $\mu=1$ and $\kappa=1$.   We will omit here the proof of existence of solutions for the system \eqref{2}  since it can  be obtained via standard finite dimensional Galerkin type approximations. For the convenience of the reader we summarize here this existence theorem.
\begin{theorem}
Let $(\ue_{0}, \te_{0}, \pe_{0})$ satisfy the conditions (ID), for some $\e>0$. Then the system \eqref{2} has a weak solution $(\ue,\te, \pe)$, with the following properties
\begin{itemize}
  \item [\bf{(i)}] $\ue\in L^{\infty}([0,T];L^{2}(\R^{3}))\cap L^{2}([0,T];\dot H^{1}(\R^{3})) $.
  \item [\bf{(ii)}] $\te\in L^{\infty}([0,T];L^{2}(\R^{3}))\cap L^{2}([0,T];\dot H^{1}(\R^{3})) $.
  \item [\bf{(iii)}] $\sqrt{\e}\pe \in L^{\infty}([0,T];L^{2}(\R^{3}))$,
\end{itemize}
for all $T>0$.
\end{theorem}

Let us now state our main result. The convergence of $\{\ue\}$ will be  described by analyzing the convergence of the associated Hodge decomposition.
\begin{theorem}
Let $(\ue, \te, \pe)$ be a sequence of weak solution in $\R^{3}$ of the system \eqref{2}, assume that the initial data satisfy (ID). Then 
\begin{itemize}
  \item [\bf{(i)}] There exists $u\in L^{\infty}([0,T];L^{2}(\R^{3}))\cap L^{2}([0,T];\dot H^{1}(\R^{3}))$ such that 
  \begin{equation*}
\ue\rightharpoonup u \quad \text{weakly in $L^{2}([0,T];\dot H^{1}(\R^{3}))$}.
\end{equation*}
  \item [\bf{(ii)}] The gradient component $Q\ue$ of the vector field $\ue$ satisfies
  \begin{equation*}
Q\ue\longrightarrow 0\quad \text{ strongly in $L^{2}([0,T];L^{p}(\R^{3}))$, for any $p\in [4,6)$}.
\end{equation*}
 \item [\bf{(iii)}] The divergence free component $P\ue$ of the vector field $\ue$ satisfies
   \begin{equation*}
P\ue\longrightarrow Pu=u\quad \text{strongly  in $L^{2}([0,T];L^{2}_{loc}(\R^{3}))$}.
\end{equation*}
\item [\bf{(iv)}] There exists $\theta\in L^{\infty}([0,T];L^{2}(\R^{3}))\cap L^{2}([0,T];\dot H^{1}(\R^{3}))$ such that 
\begin{align*}
\te&\longrightarrow \theta\quad \text{strongly  in $L^{2}([0,T];L^{2}_{loc}(\R^{3}))$,}\\
\nabla\te&\rightharpoonup \nabla\theta \quad \text{weakly  in $L^{2}([0,T]\times \R^{3})$.}
\end{align*}
 \item [\bf{(v)}] The sequence $\{\pe\}$ will converge in the sense of distribution to 
 \begin{equation*}
p=\Delta^{-1}\dive \left((u\cdot\nabla)u\right)=\Delta^{-1}tr((Du)^{2}).
\end{equation*}
\item [\bf{(vi)}] $u=Pu$ and $\theta$ are  weak solutions to the incompressible Navier Stokes Fourier system in the sense of the Definition \eqref{df}
and the following energy inequality holds for all $t\in [0,T]$,
\begin{align*}
\frac{1}{2}\int_{\R^{3}}(|u(x,t)|^{2}&+|\theta(x,t)|^{2})dx+\int_{0}^{t}\!\!\int_{\R^{3}}(|\nabla u(x,s)|^{2}+|\nabla \theta(x,s)|^{2})dxds\\&\leq 
\frac{1}{2}\int_{\R^{3}}(|u_{0}(x)|^{2}+|\theta_{0}(x)|^{2})dx.
%\label{en}
\end{align*}
\end{itemize}
\label{tM}
\end{theorem}

\section{A priori estimates}
In this section we wish to establish the a priori estimates, independent on $\e$, for the solutions of the system \eqref{2} which are necessary to prove the Theorem \ref{tM}.
 We will achieve this goal in two steps. First of all we will recover the  a priori estimates that come from the classical energy estimates related to the system \eqref{2}. Then, we get  stronger  estimates by exploiting the structure of the system. 
In fact, as we will see later on, the sequence $\pe$ satifies a  wave type equation. This will allow us to apply to $\pe$ the Strichartz estimates  \eqref{s2}, \eqref{s1} and to get in this way dispersive  bounds on $\pe$.
 \subsection{Energy estimates}
The next results concerns the  energy type estimate for the system \eqref{2}.
\begin{theorem}
Let us consider the solution $(\ue, \te, \pe)$ of the Cauchy problem for the system \eqref{2}. Assume that the hypotheses (ID) hold, then one has
\begin{equation}
\label{9}
E(t)+\int_{0}^{t}\!\!\int_{\R^{3}}(|\nabla \ue(x,s)|^{2}+|\nabla \te(x,s)|^{2})dxds=E(0),
\end{equation}
where we set
\begin{equation}
\label{8}
E(t)=\int_{\R^{3}}\left( \frac{1}{2}|\ue(x,t)|^{2}+ \frac{1}{2}|\te(x,t)|^{2}+\frac{\e}{2} |\pe(x,t)|^{2}\right)dx.
\end{equation}
\label{t1}
\end{theorem}
\begin{proof}
We multiply, as usual,  the first and second equation of the system \eqref{2} respectively by $\ue$ and $\te$ and the third by $\pe$, then we sum up and integrate by parts in space and time, hence  we get \eqref{9}.
\end{proof}
\begin{corollary}
\label{c1}
Let us consider the solution $(\ue, \te, \pe)$ of the Cauchy problem for the system \eqref{2}. Let us assume that the hypotheses (ID) hold, then  it follows
\begin{align} 
& \sqrt{\e}\pe &\quad &\text{is bounded in $L^{\infty}([0,T];L^{2}(\R^{3}))$,} \label{11}\\
& \e\pe_{t} &\quad  &\text{is relatively compact in $H^{-1}([0,T]\times \R^{3}),$}  \label{12}\\
& \nabla\ue, \nabla\te &\quad &\text{are bounded in $L^{2}([0,T]\times\R^{3}),$}  \label{13}\\
& \ue, \te&\quad&\text{are bounded in $L^{\infty}([0,T];L^{2}(\R^{3}))\cap L^{2}([0,T];L^{6}(\R^{3})),$}  \label{14}\\
&\begin{matrix}
(\ue \!\cdot\!\nabla)\ue\\ \ue\nabla\te
\end{matrix}&\quad &\text{are bounded in $L^{2}([0,T];L^{1}(\R^{3}))\!\cap \!L^{1}([0,T];L^{3/2}(\R^{3})),$}  \label{17}\\
&\begin{matrix}
(\dive \ue)\ue\\(\dive \ue)\te
\end{matrix} &\quad&\text{are bounded in $L^{2}([0,T];L^{1}(\R^{3}))\!\cap\! L^{1}([0,T];L^{3/2}(\R^{3})).$}\label{17a}
\end{align}
\end{corollary}
\begin{proof}
\eqref{11},  \eqref{12}, \eqref{13} follow from \eqref{9}, while \eqref{14} follows from \eqref{9} and Sobolev's embeddings theorems. Finally \eqref{17} and \eqref{17a} come from \eqref{13}, \eqref{14}. 
\end{proof}
\subsection{Pressure wave equation}
In this section we wish to get stronger bounds on $\pe$. First of all, we observe
that $\pe$ satisfies a wave equation. In fact, by differentiating with respect to time the equation $\eqref{2}_{3}$ and by using $\eqref{2}_{1}$ we have
$$\e\partial_{tt}\pe-\Delta \pe+\Delta\dive\ue-\dive\left((\ue\cdot\nabla)\ue+\frac{1}{2}(\dive \ue)\ue\right)=0.$$
Now, by rescaling  the time variable as
$$\tau=\frac{t}{\sqrt{\e}}$$
and consequently  the velocity, the temperature and the pressure in the following way 
\begin{equation}
\label{18}
 \ut(x,\tau)=\ue(x,\sqrt{\e}\tau), \quad  \tht(x,\tau)=\te(x,\sqrt{\e}\tau)\quad \pt(x,\tau)=\pe(x,\sqrt{\e}\tau),
\end{equation}
 we get that $\pt$ satisfies the following wave equation
\begin{equation}
\label{20 }
\partial_{\tau\tau}\pt-\Delta\pt +\Delta \dive\ut-\dive\left(\left(\ut\cdot\nabla\right)\ut +\frac{1}{2}(\dive \ut)\ut \right)=0.
\end{equation}
This structure allow us to use on $\pt$ the Strichartz estimates \eqref{s1}, \eqref{s2}.
Now the analysis for $\pt$ will follow the same line of arguments as in \cite{DM06}.
So we split $\pt$ as $\pt=\pt_{1}+\pt_{2}$ where $\pt_{1}$ and $\pt_{2}$  solve the following wave equations:
\begin{equation}
\label{21}
\begin{cases}
     \partial_{\tau\tau}\pt_{1}-\Delta\pt_{1} =-\Delta \dive\ut=F_{1} \\
     \pt_{1}(x,0)=\partial_{\tau} \pt_{1}(x,0)=0,
   \end{cases}
\end{equation}    
\begin{equation}
\label{22} 
\begin{cases}   
\displaystyle{ \partial_{\tau\tau}\pt_{2}-\Delta\pt_{2} =\dive\left(\left(\ut\cdot\nabla\right)\ut +\frac{1}{2}(\dive \ut)\ut \right)=F_{2}}\\
\pt_{2}(x,0)=\pt(x,0)\quad  \partial_{\tau}\pt_{2}(x,0)=\partial_{\tau}\pt(x,0).
 \end{cases}
 \end{equation}
Therefore we are able to prove the following theorem.
\begin{theorem}
Let us consider the solution $(\ue, \te, \pe)$ of the Cauchy problem for the system \eqref{3}. Assume that the hypotheses (ID) hold. Then we set the following estimate 
\begin{align}
\hspace{-1mm}\e^{3/8}\|\pe\|_{L^{4}_{t} W^{-2,4}_{x}}+\e^{7/8}\|\partial_{t}\pe\|_{L^{4}_{t} W^{-3,4}_{x}}&\lesssim \sqrt{\e}\|\pe_{0}\|_{L^{2}_{x}}+\|\dive\ue_{0}\|_{H^{-1}_{x}}\notag\\&+\sqrt{T}\|\dive \ue\|_{L^{2}_{t}L^{2}_{x}}\notag\\&+
\|\left(\ue\cdot\nabla\right)\ue +\frac{1}{2}(\dive \ue)\ue\|_{L^{1}_{t}L^{3/2}_{x}}.
\label{23}
\end{align}
\label{t2}
\end{theorem}
\begin{proof}
Since $\pt_{1}$ and $\pt_{2}$ are solutions of the wave equations \eqref{21}, \eqref{22}, we can apply the Strichartz estimates \eqref{s3} and \eqref{s1}, with $(x,\tau)\in \R^{3}\times\left (0,T/\sqrt \e\right)$.
Since  $\Delta ^{-1}\pt_{1}$ satisfies the equation
\begin{equation}
\label{ 24}
 \partial_{\tau\tau}(\Delta^{-1}\pt_{1})-\Delta(\Delta^{-1}\pt_{1}) =\Delta ^{-1}F_{1},
\end{equation}
then by using the Strichartz estimates \eqref{s1} we get
\begin{equation}
\|\pt_{1}\|_{L^{4}_{\tau} W^{-2,4}_{x}}+\|\partial_{\tau}\pt_{1}\|_{L^{4}_{\tau} W^{-3,4}_{x}}\lesssim \frac{\sqrt{T}}{\e^{1/4}}\|\dive \ut\|_{L^{2}_{\tau}L^{2}_{x}}.
\label{26}
\end{equation}
In the same way we have that $\Delta ^{-1/2}\pt_{2}$ satisfies the equation
\begin{equation}
\label{ 27}
 \partial_{\tau\tau}(\Delta^{-1/2}\pt_{2})-\Delta(\Delta^{-1/2}\pt_{1}) =\Delta ^{-1/2}F_{2},
\end{equation}
therefore by using  the estimate \eqref{s3} we obtain
\begin{align}
\label{29}
\|\pt_{2}\|_{L^{4}_{\tau} W^{-1,4}_{x}}+\|\partial_{\tau}\pt_{2}\|_{L^{4}_{\tau} W^{-2,4}_{x}}&\lesssim \|\pt(x,0)\|_{ H^{-1/2}_{x}}+
\|\partial_{\tau}\pt(x,0)\|_{ H^{-3/2}_{x}}\notag
\\&+\|\left(\ut\cdot\nabla\right)\ut +\frac{1}{2}(\dive \ut)\ut \|_{L^{1}_{\tau}L^{3/2}_{x},}
\end{align}
Now by using \eqref{26}, \eqref{29} it follows that $\pt$ verifies
\begin{align}
\|\pt\|_{L^{4}_{\tau} W^{-2,4}_{x}}+\|\partial_{\tau}\pt\|_{L^{4}_{\tau} W^{-3,4}_{x}}&\leq \|\pt_{1}\|_{L^{4}_{\tau} W^{-2,4}_{x}}+\|\pt_{2}\|_{L^{4}_{\tau} W^{-1,4}_{x}}\\&+
\|\partial_{\tau}\pt_{1}\|_{L^{4}_{\tau} W^{-3,4}_{x}}+\|\partial_{\tau}\pt_{2}\|_{L^{4}_{\tau} W^{-2,4}_{x}}
\notag
\\&\lesssim \|\pt(x,0)\|_{ H^{-1/2}_{x}}+
\|\partial_{\tau}\pt(x,0)\|_{ H^{-3/2}_{x}}\notag\\&
+\frac{\sqrt{T}}{\e^{1/4}}\|\dive \ut\|_{L^{2}_{\tau}L^{2}_{x},}+
\|\left(\ut\cdot\nabla\right)\ut +\frac{1}{2}(\dive \ut)\ut \|_{L^{1}_{\tau}L^{3/2}_{x}}.\notag
\label{30}
\end{align}
Finally, since
\begin{equation*}
\|\pt\|_{L^{r}((0,T/\sqrt \e );L^{q}(\R^{3}))}=\e^{-1/2r}\|\pe\|_{L^{r}([0,T];L^{q}(\R^{3}))}
\end{equation*}
 we end up with \eqref{23}.
\end{proof}
\section{Strong convergence}
In this section we conclude the proof of the  Theorem \ref{tM}. We will show the strong convergence of $Q\ue$, $P\ue$ and $\te$. In particular we will be able to prove that the gradient part of the velocity $Q\ue$ converges strongly to $0$, while the incompressible component of the velocity field $P\ue$ converges strongly to $Pu=u$, where $u$ is the limit profile as $\e\downarrow 0$ of $\ue$. Moreover we will be able to prove some time regularity of $\te$ that entails the strong compactness for $\te$ and consequently the strong convergence to $\theta$, solution of the system \eqref{1}.
\subsection{Strong convergence of $Q\ue$, $P\ue$ and $\te$}
We start this section with some easy consequences of the a priori estimates established in the previous section.
\begin{proposition}
\label{c2}
Let us consider the solution $(\ue, \te, \pe)$ of the Cauchy problem for the system \eqref{2}. Assume that the hypotheses (ID) hold. Then, as $\e\downarrow 0$, one has
\begin{align}
&\e\pe\longrightarrow 0 &\quad& \text{strongly in $L^{\infty}([0,T];L^{2}(\R^{3}))\cap L^{4}([0,T];W^{-2,4}(\R^{3}))$,}\label{31}\\
&\dive \ue \longrightarrow 0 &\quad& \text{strongly in $ W^{-1,\infty}([0,T];L^{2}(\R^{3}))\cap L^{4}([0,T];W^{-3,4}(\R^{3}))$}.\label{32}
\end{align}
\end{proposition}
\begin{proof}
\eqref{31}, \eqref{32} follow  from the estimates \eqref{11}, \eqref{23} and the second equation of the system \eqref{3}. 
\end{proof}

Now, we wish to show that the gradient part of the velocity field $Q\ue$ goes strongly to $0$ as $\e\downarrow 0$. As we will see in the next proposition, this will be a consequence of the estimate \eqref{23} and of the Lemma \ref{ly}.
\begin{proposition}
Let us consider the solution $(\ue, \te, \pe)$ of the Cauchy problem for the system \eqref{2}. Assume that the hypotheses (ID) hold. Then  as $\e\downarrow 0$,
\begin{equation}
Q\ue \longrightarrow 0 \quad \text{strongly in $ L^{2}([0,T];L^{p}(\R^{3}))$ for any $p\in [4,6)$ }.
\label{33}
\end{equation}
\label{p2}
\end{proposition}
\begin{proof}
In order to prove the Proposition \ref{p2} we split $Q\ue$ as follows
\begin{equation*}
\|Q\ue\|_{L^{2}_{t}L^{p}_{x}}\leq \|Q\ue-Q\ue\ast j_{\alpha}\|_{L^{2}_{t}L^{p}_{x}}+\|Q\ue\ast j_{\alpha}\|_{L^{2}_{t}L^{p}_{x}}=J_{1}+J_{2},
\end{equation*}
where $j_{\alpha}$ is the smoothing kernel defined in Lemma \ref{ly}.
Now we estimate separately $J_{1}$ and $J_{2}$. For $J_{1}$ by using \eqref{y1} we get
\begin{equation}
\label{50}
J_{1}\leq \alpha^{1-3\left(\frac{1}{2}-\frac{1}{p}\right)}\left(\int_{0}^{T}\|\nabla Q\ue(t)\|_{L^{2}_{x}}^{2} dt\right)^{1/2}\leq \alpha^{1-3\left(\frac{1}{2}-\frac{1}{p}\right)}\|\nabla \ue\|_{L^{2}_{t}L^{2}_{x}}.
\end{equation}
Hence from the identity $Q\ue=-\e^{1/8}\nabla\Delta^{-1}\e^{7/8}\partial_{t}p$ and by the inequality \eqref{y2} we get  $J_{2}$ satisfies the following estimate
\begin{equation}
J_{2}\leq \e^{1/8}\alpha^{-2-3\left(\frac{1}{4}-\frac{1}{p}\right)}T^{1/4}\|\e^{7/8}\partial_{t}p\|_{L^{4}_{t}W^{-3,4}_{x}}.
\label{51}
\end{equation}
Therefore,  summing up \eqref{50} and \eqref{51} and by using \eqref{13} and \eqref{23}, we conclude  for any $p\in [4,6)$ that
\begin{equation}
\|Q\ue\|_{L^{2}_{t}L^{p}_{x}}\leq C\alpha^{1-3\left(\frac{1}{2}-\frac{1}{p}\right)}+C_{T}\e^{1/8}\alpha^{-2-3\left(\frac{1}{4}-\frac{1}{p}\right)}.
\end{equation}
Finally we choose $\alpha$ in terms of $\e$ in order that the two terms in the right hand side of the previous inequality have the same order, namely
\begin{equation}
\alpha=\e^{1/18}.
\end{equation}
Therefore we obtain
\begin{equation*}
\displaystyle{\|Q\ue\|_{L^{2}_{t}L^{p}_{x}}\leq C_{T}\e^{ \frac{6-p}{36p}}\quad \text{for any $p\in [4,6)$.}}
\end{equation*}
\end{proof}

It remains to prove the strong compactness of the incompressible component of the velocity field and of the temperature.  To achieve this goal we need to recall here, the following theorem (see \cite{Si}).
\begin{theorem}
Let be $\mathcal{F}\subset L^{p}([0,T];B)$,  $1\leq p<\infty$, $B$ a Banach space. $\mathcal{F}$ is relatively compact in  $L^{p}([0,T];B)$ for $1\leq p<\infty$, or in $C([0,T];B)$ for $p=\infty$ if and only if 
\begin{itemize}
\item[{\bf (i)}]
$\displaystyle{\left\{\int_{t_{1}}^{t_{2}}f(t)dt,\ f\in B\right\}}$ is relatively compact in $B$, $0<t_{1}<t_{2}<T$,
\item[{\bf (ii)}]
$\displaystyle{\lim_{h\to 0}\|f(t+h) - f(t)\|_{L^{p}([0, T-h];B)}=0}$ uniformly for any $f \in \mathcal{F}$.
\end{itemize}
\label{S}
\end{theorem}

The compactness can be obtained by looking at some time regularity properties of $P\ue$ and $\te$ and by using the Theorem \ref{S}, but before we need to prove  the following lemmas.
\begin{lemma}
\label{l2}
Let us consider the solution $(\ue, \te, \pe)$ of the Cauchy problem for the system \eqref{2}. Assume that the hypotheses (ID) hold. Then for all $h\in(0,1)$, we have
\begin{equation}
\label{34t}
\|\te(t+h)-\te(t)\|_{L^{2}([0,T]\times \R^{3})}\leq C_{T}h^{1/5}.
\end{equation}
\end{lemma}
\begin{proof}
Let us set $w^{\e}=\te(t+h)-\te(t)$, we have
\begin{align}
\hspace{-0,25 cm}\|\te(t+h)-\te(t)\|^{2}_{L^{2}([0,T]\times \R^{3})}&=\int_{0}^{T}\!\!\int_{\R^{3}}dtdx w^{\e}\cdot(w^{\e}-w^{\e}\ast j_{\alpha})\notag\\&+\int_{0}^{T}\!\!\int_{\R^{3}}dtdx w^{\e}\cdot(w^{\e}\ast j_{\alpha})=I_{1}+I_{2}.
\label{35t}
\end{align}
By using \eqref{y1} we can estimate $I_{1}$ in the following way
\begin{align}
I_{1}&\leq \|w^{\e}\|_{L^{\infty}_{t}L^{2}_{x}}\int_{0}^{T}\|w^{\e}(t)-(w^{\e}\ast j_{\alpha})(t)\|_{L^{2}_{x}}dt\notag\\&\lesssim \alpha T^{1/2}\|\te\|_{L^{\infty}_{t}L^{2}_{x}}\|\nabla\te\|_{L^{2}_{t,x}}.
\label{36t}
\end{align}
Let us reformulate $w^{\e}$ in integral form by using the equation $\eqref{2}_{2}$, hence
\begin{align}
\hspace{-0.3cm}I_{2}\leq\left|\int_{0}^{T}\!\!\!dt\!\!\int_{\R^{3}}\!\!\!dx \!\!\int_{t}^{t+h}\!\!\!\!\!ds(\Delta \te-\ue\cdot\nabla\te -\frac{1}{2}\te(\dive \ue)(s,x)\cdot (w^{\e}\ast j_{\alpha})(t,x)\right|.
\label{37t}
\end{align}
Then integrating by parts and by using \eqref{y2}, with $s=0$, $p=\infty$ and $q=2$, we deduce
\begin{align}
I_{2}&\leq h\|\nabla\te\|^{2}_{L^{2}_{t,x}}+C\alpha^{-3/2}T^{1/2}\|\te\|_{L^{\infty}_{t}L^{2}_{x}}\left(\!h\!\int_{t}^{t+h}\!\!\!\|\ue\cdot\nabla\te -\frac{1}{2}(\dive \ue)\te\|^{2}_{L^{1}_{x}}ds\!\!\right)^{1/2}\notag\\
&\leq h\|\nabla\te\|^{2}_{L^{2}_{t,x}}\!\!+C\alpha^{-3/2}T^{1/2}h\|\ue\|_{L^{\infty}_{t}L^{2}_{x}}\|\ue\!\cdot\!\nabla\te \!-\frac{1}{2}(\dive \ue)\te\|_{L^{2}_{t}L^{1}_{x}}.
\label{38t}
\end{align}
Summing up $I_{1}$, $I_{2}$ and by taking into account \eqref{13}, \eqref{14}, \eqref{17}, \eqref{17a}, we have
\begin{equation}
\|\te(t+h)-\te(t)\|^{2}_{L^{2}([0,T]\times \R^{3})}\leq C(\alpha T^{1/2}+h\alpha^{-3/2}T^{1/2}+h),
\label{39t}
\end{equation}
by choosing $\alpha=h^{2/5}$, we end up with \eqref{34t}.
\end{proof}

As for $\te$ we prove the equicontinuity in time for the sequence $P\ue$.
\begin{lemma}
\label{l1}
Let us consider the solution $(\ue, \te, \pe)$ of the Cauchy problem for the system \eqref{2}. Assume that the hypotheses (ID) hold. Then for all $h\in(0,1)$, we have
\begin{equation}
\label{34}
\|P\ue(t+h)-P\ue(t)\|_{L^{2}([0,T]\times \R^{3})}\leq C_{T}h^{1/5}.
\end{equation}
\end{lemma}
\begin{proof}
Let us set $z^{\e}=\ue(t+h)-\ue(t)$, we have
\begin{align}
\hspace{-0,25 cm}\|P\ue(t+h)-P\ue(t)\|^{2}_{L^{2}([0,T]\times \R^{3})}&=\int_{0}^{T}\!\!\int_{\R^{3}}dtdx(Pz^{\e})\cdot(Pz^{\e}-Pz^{\e}\ast j_{\alpha})\notag\\&+\int_{0}^{T}\!\!\int_{\R^{3}}dtdx(Pz^{\e})\cdot(Pz^{\e}\ast j_{\alpha})=I_{1}+I_{2}.
\label{35}
\end{align}
By using \eqref{y1} we can estimate $I_{1}$ in the following way
\begin{equation}
I_{1}\lesssim \alpha T^{1/2}\|\ue\|_{L^{\infty}_{t}L^{2}_{x}}\|\nabla\ue\|_{L^{2}_{t,x}}.
\label{36}
\end{equation}
Let us reformulate $Pz^{\e}$ in integral form by using the equation $\eqref{2}_{1}$, hence
\begin{align}
\hspace{-0.3cm}I_{2}\leq\left|\int_{0}^{T}\!\!\!dt\!\!\int_{\R^{3}}\!\!\!dx \!\!\int_{t}^{t+h}\!\!\!ds(\Delta \ue-\left(\ue\cdot\nabla\right)\ue -\frac{1}{2}\ue(\dive \ue)(s,x)\cdot (Pz^{\e}\ast j_{\alpha})(t,x)\right|.
\label{37}
\end{align}
Then integrating by parts and by using \eqref{y2}, with $s=0$, $p=\infty$ and $q=2$, we deduce
\begin{equation}
I_{2}\leq h\|\nabla\ue\|^{2}_{L^{2}_{t,x}}+C\alpha^{-3/2}T^{1/2}h\|\ue\|_{L^{\infty}_{t}L^{2}_{x}}\|\left(\ue\cdot\nabla\right)\ue -\frac{1}{2}(\dive \ue)\ue\|_{L^{2}_{t}L^{1}_{x}}.
\label{38}
\end{equation}
Summing up $I_{1}$, $I_{2}$ and by taking into account \eqref{13}, \eqref{14}, \eqref{17}, \eqref{17a}, we have
\begin{equation}
\|P\ue(t+h)-P\ue(t)\|^{2}_{L^{2}([0,T]\times \R^{3})}\leq C(\alpha T^{1/2}+h\alpha^{-3/2}T^{1/2}+h),
\label{39}
\end{equation}
by choosing $\alpha=h^{2/5}$, we end up with \eqref{34}.
\end{proof}

\begin{corollary}
Let us consider the solution $(\ue, \te, \pe)$ of the Cauchy problem for the system \eqref{2}. Assume that the hypotheses (ID) hold. Then  as $\e\downarrow 0$
\label{c3}
\begin{align}
P\ue &\longrightarrow Pu, &\qquad& \text{strongly in $L^{2}(0,T;L^{2}_{loc}(\R^{3}))$}.
\label{43}\\
\te &\longrightarrow \theta, &\qquad& \text{strongly in $L^{2}(0,T;L^{2}_{loc}(\R^{3}))$}.
\label{43t}\\
\nabla\te &\rightharpoonup \nabla\theta, &\qquad& \text{weakly in $L^{2}((0,T)\times \R^{3})$}.
\label{43tt}
\end{align}
\end{corollary}
\begin{proof}
By using the Lemma \ref{l1}, the Lemma \ref{l2} and the Theorem \ref{S} and the Proposition \ref{p2} we get \eqref{43}, \eqref{43t}, while \eqref{43tt} is a consequence of \eqref{13}.
\end{proof}
\subsection{Proof of the Theorem \ref{tM}}
\begin{itemize}
\item[{\bf (i)}] It follows from the estimate \eqref{14}.
\item[{\bf (ii)}] It is a consequence of  the Proposition \ref{p2}.
\item[{\bf (iii)}] By taking into account the decomposition $\ue=P\ue+Q\ue$, by the Corollary \ref{c3} and the Proposition \ref{p2} we have that
\begin{equation*}
P\ue\longrightarrow u \qquad \text{strongly in $L^{2}([0,T];L^{2}_{loc}(\R^{3}))$.}
\end{equation*}
\item[{\bf (iv)}] Follows from the estimate \eqref{14} and by using the Corollary \ref{c3}
\item[{\bf (v)}]  Let us apply the Leray projector $Q$ to the equation $\eqref{3}_{1}$, then it follows
\begin{equation}
\label{54}
\nabla \pe =\Delta Q\ue- Q\left(\dive(\ue\otimes\ue) +\frac{3}{2} \ue \dive Q\ue\right).
\end{equation}
Now by choosing a test function $\varphi \in H^{1}_{t}W^{2,4/3}_{x}\cap C^{0}_{t}W^{1,4}_{x}\cap L^{2}_{t}H^{1}_{x}$ and by taking into account \eqref{13}, \eqref{33}, \eqref{43}, we get, as $\e \downarrow 0$, 
\begin{align}
\langle\ue \dive Q\ue, Q\varphi\rangle&\leq\|Q\ue\|_{L^{2}_{t}L^{4}_{x}} \|\nabla\ue\|_{L^{2}_{t}L^{2}_{x}}\|Q\varphi\|_{L^{\infty}_{t}L^{4}_{x}} \notag\\&+ \|Q\ue\|_{L^{2}_{t}L^{4}_{x}} \|\ue\|_{L^{\infty}_{t}L^{2}_{x}}\|\nabla Q\varphi\|_{L^{2}_{t}L^{4}_{x}} \rightarrow 0, 
\end{align}
\begin{align}
\langle \dive(\ue\otimes\ue),Q\varphi\rangle&=\langle \dive(P\ue\otimes P\ue),Q\varphi\rangle+\langle \dive(Q\ue\otimes Q\ue),Q\varphi\rangle\notag\\&+\langle \dive(P\ue\otimes Q\ue),Q\varphi\rangle\notag+\langle \dive(Q\ue\otimes Q\ue),Q\varphi \rangle\notag\\&\rightarrow \langle \dive(Pu \otimes Pu),Q\varphi\rangle=\langle Q\dive((Pu\cdot\nabla)Pu) , \varphi\rangle.
\end{align}
So  as $\e \downarrow 0$ one has, 
\begin{equation}
\label{40}
\langle \nabla \pe , \varphi \rangle \longrightarrow\langle \nabla\Delta^{-1}\dive((u\cdot\nabla)u) , \varphi\rangle. 
\end{equation}
\item[{\bf(vi)}] Before passing into the limit in the system \eqref{2} the convergence of the nonlinear term  $\ue\cdot\nabla\te$ deserves a little discussion. By using again for $\ue$ the associated Hodge decomposition, namely $\ue=P\ue+Q\ue$ and by taking into account \eqref{33}, \eqref{43}, \eqref{43t}, \eqref{43tt} one can prove that, for any $p\in[4,6)$,
\begin{align}
Q\ue\cdot\nabla\te &\longrightarrow 0 &\qquad&\text{strongly in $L^{1}([0,T];L^{\frac{2p}{p+2}}(\R^{3}))$,}\\
P\ue\cdot\nabla\te &\longrightarrow u\cdot\nabla\theta &\qquad&\text{strongly in $L^{1}([0,T]\times \R^{3})$}. 
\end{align}

Now we can pass  into the limit inside the system \eqref{2} and we get  $u$ and $\theta$ satisfy the following equation in $\mathcal{D}'([0,T]\times \R^{3})$
\begin{align}
P(\partial_{t} u-\Delta u+(u\cdot\nabla)u)=0,\label{55}\\
\partial_{t} \theta-\Delta \theta+u\cdot\nabla\theta=0.\label{55t} 
\end{align}
Finally we prove the energy inequality. By using the weak lower semicontinuity of the weak limits, the hypotheses (ID) and denoting 
by $\chi$ the weak-limit of $\sqrt{\e}\pe $,  we have
\begin{align}
&\int_{\R^{3}}\frac{1}{2}\left(|\chi(x,t)|^{2}+|u(x,t)|^{2}+|\theta(x,t)|^{2}\right)dx\notag\\&+\int_{0}^{t}\!\!\int_{\R^{3}}\left(|\nabla u(x,s)|^{2}+|\nabla \theta(x,s)|^{2}\right)dxds\notag\\&\leq
\liminf_{\e\to 0}\left(\frac{1}{2}\int_{\R^{3}}(\e|\pe(x,t)|^{2}+|\ue(x,t)|^{2}+|\te(x,t)|^{2})dx\right)\notag\\
&+\liminf_{\e\to 0}\left(\int_{0}^{t}\!\!\int_{\R^{3}}(|\nabla \ue(x,s)|^{2}+|\nabla \te(x,s)|^{2})dxds\right)\notag\\&=\liminf_{\e\to 0}\int_{\R^{3}}\frac{1}{2}\left(\e|\pe_{0}(x)|^{2}+|\ue_{0}(x)|^{2}+|\te_{0}(x)|^{2}\right)dx\notag\\
&=\int_{\R^{3}}\frac{1}{2}(|u_{0}(x)|^{2}+|\theta_{0}(x)|^{2})dx.
\end{align}
\end{itemize}

\end{document}